# FINITARY ISOMORPHISMS OF BROWNIAN MOTIONS

BY ZEMER KOSLOFF[1] AND TERRY SOO[2]

[1]*Einstein Institute of Mathematics, Hebrew University of Jerusalem, zemer.kosloff@mail.huji.ac.il*
[2]*Department of Statistical Science, University College London, math@terrysoo.com*

Ornstein and Shields (*Advances in Math.* **10** (1973) 143–146) proved that Brownian motion reflected on a bounded region is an infinite entropy Bernoulli flow, and, thus, Ornstein theory yielded the existence of a measure-preserving isomorphism between any two such Brownian motions. For fixed $h > 0$, we construct by elementary methods, isomorphisms with almost surely finite coding windows between Brownian motions reflected on the intervals $[0, qh]$ for all positive rationals $q$.

**1. Introduction.** A measure-preserving flow $(K, \kappa, T)$ is *Bernoulli* if for each time $t > 0$, the discrete-time system $(K, \kappa, (T_{nt})_{n \in \mathbb{Z}})$ is isomorphic to a Bernoulli shift; see Section 2.4 for background on Bernoulli shifts. The entropy of a flow is given by the Kolmogorov–Sinai entropy of its time-one map. The remarkable theory developed by Ornstein and his collaborators tells us that all infinite entropy Bernoulli flows are isomorphic [22]. Ornstein and Shields [23] proved that mixing Markov shifts of kernel type are Bernoulli, which include Brownian motion on bounded reflecting regions, and also irreducible continuous-time Markov chains which were earlier shown to be Bernoulli by Feldman and Smorodinsky [6].

The isomorphisms provided by Ornstein theory are given by an abstract existence proof that offer no control over the coding windows. When the coding window is finite, almost surely, then we say that the isomorphism is finitary; we will define this concept more carefully later in Section 2.4. Recently, Soo and Wilkens [32] constructed finitary isomorphisms between any two Poisson point processes, and Soo [31] constructed finitary isomorphisms between Poisson point process and a restricted class of continuous-time finite state Markov chains. We will construct by elementary methods finitary isomorphisms between reflected Brownian motions.

In what follows we consider reflected and periodic Brownian motions arising from Brownian motion with variance $\sigma^2 = 1$. We will discuss the general case $\sigma > 0$ in Section 6.2.

THEOREM 1. *Let $h > 0$. Let $q$ be a positive rational. A stationary Brownian motion reflected on the interval $[0, qh]$ is finitarily isomorphic to a stationary Brownian motion reflected on $[0, h]$; moreover, the isomorphism also has a finitary inverse.*

In order to prove Theorem 1, we will use the closely related model of a stationary periodic Brownian motion on a bounded region as an intermediate process.

THEOREM 2. *Let $h > 0$. A stationary periodic Brownian motion on $[0, 2h]$ is finitarily isomorphic to a stationary Brownian motion reflected on $[0, h]$; moreover, the isomorphism also has a finitary inverse.*

THEOREM 3. *Let $h > 0$ and $q$ be a positive rational. A stationary periodic Brownian motion on the interval $[0, qh]$ is finitarily isomorphic to a stationary periodic Brownian motion on $[0, h]$; moreover, the isomorphism also has a finitary inverse.*

PROOF OF THEOREM 1. Immediate from Theorems 2 and 3. □

We note that standard scaling properties enjoyed by Brownian motions or Poisson point processes cannot be exploited in an obvious way to produce isomorphisms since such an exploitation will violate equivariance; see Section 2.1. The proof of Theorem 2 relies on the key observation that a renewal point process, corresponding to certain hitting times, appears as a factor in *both* reflected and periodic Brownian motions; furthermore, the laws of these Brownian motions conditioned on this renewal point process has a nice product structure which can be encoded by independent and identically distributed (i.i.d.) random variables that are uniformly distributed on $[0, 1]$. The renewal point process plays the role of *markers* and the conditioned Brownian paths *fillers* in the sense of Keane and Smorodinsky [14]; see also Section 2.4.

It does not appear that our proof of Theorem 3 (and thus Theorem 1) extends in an obvious way that would allow one to remove the condition that the scaling $q$ is rational.

QUESTION 1. *Can the restriction that the scaling $q > 0$ is rational be removed in Theorems 1 and 3?*

Poisson point processes can be considered to be canonical examples of infinite entropy Bernoulli flows. We ask the following question:

QUESTION 2. *Do there exist finitary isomorphisms between Brownian motions reflected on bounded regions and Poisson point processes?*

The tools we develop to prove Theorem 1 easily yield the following partial answer to Question 2:

THEOREM 4. *Poisson point processes are finitary factors of Brownian motions reflected on bounded regions.*

## 2. Preliminary definitions.

2.1. *Factors and isomorphisms.* Let $(K, \kappa, T)$ and $(K', \kappa', T')$ be two measure-preserving flows. We say that $\phi : K \to K'$ is *factor* from $(K, \kappa, T)$ to $(K', \kappa', T')$ if on a set of full measure we have $\phi \circ T_s = T'_s \circ \phi$ for all $s \in \mathbb{R}$ (*equivariance*) and $\mu' = \mu \circ \phi^{-1}$; if $\phi$ is injective, then $\phi^{-1} : K' \to K$ is a factor $(K, \kappa, T)$ to $(K', \kappa', T')$, the flows are *isomorphic* and $\phi$ is an *isomorphism*. In what follows we will define precisely how Poisson point processes, reflected Brownian motions, and related processes from probability theory give rise to measure-preserving flows.

2.2. *Simple stationary point processes and renewal processes.* A simple point process on $\mathbb{R}$ is a random variable $\Pi$ that takes values on the space of all Borel simple point measures denoted by $\mathbb{M}$. Thus, if $\mu \in \mathbb{M}$, then it is a countable sum of Dirac point masses where multiplicity is not allowed. The point processes we consider will only have isolated points, almost surely. Let $T = (T_t)_{t \in \mathbb{R}}$ be the group of translations of $\mathbb{R}$ which act on $\mathbb{M}$ via $T_t(\mu) = \mu \circ T_{-t}$. The point process $\Pi$ is *stationary* if the system associated with $\Pi$ given by $(\mathbb{M}, \mathbb{P}(\Pi \in \cdot), T)$ is measure preserving. In what follows it will be convenient to identify $\mathbb{M}$ as the collection of countable subsets of $\mathbb{R}$. Using the ordering of $\mathbb{R}$, a realization of $\Pi$ will be written as $[\Pi(\omega)] = \{p_i(\omega)\}_{i \in \mathbb{Z}}$ where $p_0(\omega) \leq 0$ is the closest point to 0 from the left and $p_i(\omega) < p_j(\omega)$ for all $i < j$.

A stationary point process $\Pi$ is a *renewal point process* if the random variables $(p_{i+1} - p_i)_{i \in \mathbb{Z}^+}$ form an i.i.d. sequence; in this case the law of $p_2 - p_1$ is called *the lifetime distribution* of $\Pi$. Note that $p_1 - p_0$ does *not* have the lifetime distribution due to *size biasing*. The most notable example of a renewal process is the Poisson point process on $\mathbb{R}$ with intensity $\lambda > 0$ which is a renewal point process whose lifetime distribution is an exponential distribution of rate $\lambda$. If the lifetime distribution of the renewal point process is nonatomic, then we will say that $\Pi$ is a *renewal point process with a continuous lifetime distribution*; all the renewal point processes considered in this paper will be of this type.

We will also consider marked point processes, where the markings are independent of the points but may have certain dependencies within themselves. If $\mu \in \mathbb{M}$, we write $[\mu] = \{p_i\}_{i \in \mathbb{Z}} \subset \mathbb{R}$ to be the points of $\mu$. Let $A$ be a set; we will mainly be concerned with the case that $A$ is finite or the unit interval $[0, 1]$. We will refer to elements of $A$ as *marks* or *colours*. Let $(\mu, x) \in \mathbb{M} \times A^{\mathbb{Z}}$; we will think of the points of $\mu$ as marked by the corresponding term in the sequence $x$ so that if $[\mu] = \{p_i\}_{i \in \mathbb{Z}}$, then

$$[(\mu, x)] = \{(p_i, x_i)\}_{i \in \mathbb{Z}}$$

is the sequence of *marked points*. Sometimes, we will say that the point $p_i$ *receives* the colour $x_i$. We will define translations of $\mathbb{M} \times A^{\mathbb{Z}}$ which preserve the markings. Let $\theta$ be a translation of $\mathbb{R}$. We set $\theta(\mu, x) = (\theta\mu, \sigma^n x)$, where $\sigma$ is the usual left shift and $n \in \mathbb{Z}$ is such that $\theta(p_n) \leq 0$ is the point of $\theta\mu$ that is closest to the origin from the left.

We say that $(\Gamma, \mathfrak{X})$ is a *marked point process* if $\Gamma$ is a point process and $\mathfrak{X}$ is a sequence of corresponding marks. We say that $\mathfrak{X}$ is the *skeleton*, and we say that the point process is *independently marked* if $\mathfrak{X}$ is independent of $\Gamma$. All the marked point processes we consider will be independently marked. Note that a stationary irreducible continuous-time finite state Markov chain where the states have the same holding rate has a presentation as an independently marked Poisson point process where the distribution of the marks is a discrete-time Markov chain.

2.3. *Brownian motions.* The Brownian motions we consider will be random variables taking values on the space $C(\mathbb{R}, I)$ of all continuous functions from $\mathbb{R}$ to an interval $I$ which will either be all of $\mathbb{R}$ or bounded. Let $T = (T_t)_{t \in \mathbb{R}}$ be the group of translations of $\mathbb{R}$ which act on $C(\mathbb{R}, I)$ via $T_t(f) = f \circ T_{-t}$. Let $W = (W_t)_{t \in \mathbb{R}}$ be a stationary Brownian motion reflected on the bounded interval $I$. The *reflected Brownian* measure-preserving system associated with $W$ is given by $(C(\mathbb{R}, I), \mathbb{P}(W \in \cdot), T)$. Similarly, if $Z = (Z_t)_{t \in \mathbb{R}}$ is a periodic Brownian motion on the bounded interval $I$, then the *periodic Brownian* measure-preserving system associated with $Z$ is given by $(C(\mathbb{R}, I), \mathbb{P}(Z \in \cdot), T)$. We will give precise definitions for $W$ and $Z$ in Sections 4 and 5.

2.4. *Finitary factors.* Recall that a *Bernoulli shift* is the measure-preserving system given by $(A^{\mathbb{Z}}, p^{\mathbb{Z}}, \sigma)$, where $A$ is a finite set, $p$ is a probability measure on $A$ and $\sigma$ is the usual left shift. The Kolmogorov–Sinai entropy of a Bernoulli shift is given by the Shannon entropy $H(p)$ and is an isomorphism invariant [18, 29]. Sinai [30] proved that a Bernoulli shift is a factor of another Bernoulli shift of greater or equal entropy. In the debut application of his theory, Ornstein [21] proved that two Bernoulli shifts of the same entropy are isomorphic. Sinai's and Ornstein's constructions are different and do not produce finitary isomorphisms.

In a landmark paper of Keane and Smorodinsky [15], they constructed, by a marker-filler method, an explicit isomorphism $\phi: A^{\mathbb{Z}} \to B^{\mathbb{Z}}$ between any two Bernoulli shifts $(A^{\mathbb{Z}}, p^{\mathbb{Z}}, \sigma)$ and $(B^{\mathbb{Z}}, q^{\mathbb{Z}}, \sigma)$ with $H(p) = H(q)$; moreover, $\phi$ is *finitary* so that there exists a measurable function $w : A^{\mathbb{Z}} \to \mathbb{N} \cup \{\infty\}$ such that $w$ is $p^{\mathbb{Z}}$-almost surely finite and $\phi(x)_0 = \phi(x')_0$ if $x|_{(-w(x), w(x))} = x'|_{(-w(x), w(x))}$, where $x|_I \in A^I$ denotes the restriction of $x$ to the coordinates in $I \subset \mathbb{Z}$. Sometimes, we will refer to the function $w$ as a *coding window*. Thus, if $\phi$ is finitary, then the coding window of $\phi$ is finite, almost surely. The map $\phi$ constructed by Keane and Smorodinsky also has a finitary inverse.

It is interesting to note the first example of a nontrivial isomorphism between two Bernoulli shifts is due to Mešalkin [20] and his isomorphism is finitary with a finitary inverse; see also [17], page 212. His beautiful and elementary example continues to have a lasting influence in a number of areas of ergodic theory and probability theory, including applications to the isomorphism problem [8, 12], shift-coupling, palm theory, matching [11] and generation of random stationary graphs [2].

Typically, one has no control over whether the coding windows of the factors and isomorphisms coming from Ornstein theory will be finite. Keane and Smorodinsky's work [14–16] inspired the hope that a theory of finitary isomorphisms could be developed in parallel to Ornstein theory [27, 28]. However, many of the celebrated consequences of Ornstein theory do not yet have finitary counterparts. Moreover, Smorodinsky exhibited an example of a countable state Markov chain with finite entropy that is Bernoulli but is not finitarily isomorphic to a Bernoulli shift [19]. In this paper, we continue to provide finitary counterparts in the infinite entropy setting.

Let $I \subset \mathbb{R}$. For $\mu \in \mathbb{M}$, we let $\mu|_I(\cdot) := \mu(\cdot \cap I)$ denote the restriction of $\mu$ to $I \subset \mathbb{R}$, and we let $f|_I \in \mathbb{R}^I$ denote the restriction of $f$ to $I$. We say that two elements of the above spaces *agree* on $I$ if their restrictions agree on $I$. Thus, we extend Keane and Smorodinsky's notion of a finitary factor in the obvious way. Let $\Pi$ be a Poisson point process, and consider a factor map $\phi$ taking $\Pi$ to a Brownian motion reflected on $[a, b]$. We say that a factor map $\phi : \mathbb{M} \to C(\mathbb{R}, [a, b])$ is *finitary* if there exists a measurable function $w : \mathbb{M} \to \mathbb{N} \cup \{\infty\}$ such that $w(\Pi) < \infty$ almost surely and $\phi(\mu)|_{(0,1)} = \phi(\mu')|_{(0,1)}$ if $\mu$ and $\mu'$ agree on $(-w(\mu), w(\mu))$. We similarly define *finitary* factors in other contexts of interest. In particular, in the case of marked point processes we say that $(\mu, x), (\mu', x') \in \mathbb{M} \times A^{\mathbb{Z}}$ *agree* on $I \subset \mathbb{R}$ if the point measures $\mu$ and $\mu'$ agree on $I$ and the points in $I$ receive the same marks.

## 3. Finitary isomorphisms between renewal processes with continuous lifetime distribution and their marked versions.
In this section we will prove some variations of results proved in [31] for Poisson point processes. We will only consider finite state Markov chains in this paper.

LEMMA 5. *Consider an independently marked stationary renewal point process with a continuous lifetime distribution and a skeleton that is an irreducible finite state Markov chain with at least two states. There exist a finitary isomorphism from this marked point process to a marked point processes, with the same law, where the points, additionally, independently receive i.i.d. marks that are uniformly distributed in* [0, 1]; *moreover, the isomorphism has a finitary inverse.*

REMARK 6. We will on occasion consider marked point processes that receive additional marks as in Lemma 5. In this context, in order to avoid confusion in proofs, we will often refer to the original marks as colours.

PROOF OF LEMMA 5. Let $(\Gamma, \mathfrak{X})$ be the coloured point process with colours $S$, so that the Markov chain $\mathfrak{X}$ takes values in the finite state space $S$. Let $r \in S$. Since the chain is irreducible, there exists states $r_1, \ldots, r_n \in S \setminus \{r\}$, such that $n \geq 1$ is minimal and with nonzero probability if the Markov chain is in state $r$ it returns to state $r$ in $n+1$ transitions via the path $\rho := (r, r_1, \ldots, r_n, r)$. We call the interval $[a, b]$ a *rho-interval* if $[\Gamma] \cap [a, b]$ is given by $\{a, t_1, \ldots, t_n, b\}$ where $a < t_1 < \cdots < t_n < b$ and the points are coloured correspondingly so that the points $a$ and $b$ receive the colour $r$ and the point $t_i$ receives the colour $r_i$. We say that a point $t$ of $\Gamma$ *reports* to the rho-interval $[a, b]$ if $b \leq t$ and $b$ is largest such right endpoint with this property; we also order all the points that report to the interval in increasing order.

A *rho-distribution on* $[a, b]$ is given by the law of $\Gamma|_{[a,b]}$, conditioned that the interval $[a, b]$ is a rho-interval. Thus, a realization of a rho-distribution on $[a, b]$ will always contain the points $a$ and $b$ and exactly $n \geq 1$ points in between $a$ and $b$. Since $(\Gamma, \mathfrak{X})$ is an independently coloured renewal process, we have that, conditional on the locations of rho-intervals of $(\Gamma, \mathfrak{X})$, the point process $\Gamma'$, formed by resampling independently on each rho-interval, according to the corresponding rho-distribution, has the same law as $\Gamma$. We refer to the points of $\Gamma$ as the *source* points for $\Gamma'$ so that the source points are resampled in the rho-intervals. If $\mathfrak{X}'$ is a skeleton such that the source points of $\Gamma'$ receive the same colour as they did under $\mathfrak{X}$ and the resampled points receive the same colour as their corresponding source points under $\mathfrak{X}$, then $(\Gamma, \mathfrak{X}) \stackrel{d}{=} (\Gamma', \mathfrak{X}')$. Note that $\mathfrak{X}$ may not be the same as $\mathfrak{X}'$ as the sequence of colours may be shifted, although the ordering is preserved.

By the Borel isomorphism theorem for probability spaces [33], Theorem 3.4.23, and the assumption that the lifetime distribution is continuous, for each integer $k \geq 1$ and every interval $[a, b]$ it follows that there exists a bijective mapping $\phi$ such that

$$\phi(\Xi) \stackrel{d}{=} (\Xi, U),$$

where $\Xi$ has the rho-distribution on $[a, b]$ and $U = (U_1, \ldots, U_k)$ are i.i.d. and uniformly distributed on $[0, 1]$ and independent of $\Xi$.

Hence, given all the rho-intervals of $(\Gamma, \mathfrak{X})$, by resampling on a rho-interval with $k$ points that report to it, using the mapping $\phi$, we additionally obtain $k$ independent random variables that are uniformly distributed on $[0, 1]$ that can be assigned its $k$ points. By independently resampling on every $\rho$ interval in this manner, we obtain the desired result. □

REMARK 7. Notice that in proof of Lemma 5, we make use of the original marks in our construction. A version of Lemma 5 is proved in [31], Proposition 4, for *unmarked* Poisson point process, and the proof requires more involved machinery from Holroyd, Lyons and Soo [9], developed for Poisson point processes.

If $\Pi$ is a stationary point process, then the *alternating* point process corresponding to $\Pi$ is obtained by colouring the points of $\Pi$ red and blue successively, where the *two* such *global* choices occur with equal probability independently of $\Pi$. Alternating point process will play a key role in our proof of Theorem 2. It is not known whether alternating Poisson point processes can be obtained as finitary factors of Poisson point processes [31], Question 1; in addition, as a factor of a Poisson point process it is impossible to colour these same points red and blue to obtain an alternating Poisson point process [10], Lemma 11. Ornstein theory, via Feldman and Smorodinsky [6], assures us that alternating Poisson point processes are Bernoulli, despite the apparent inherent periodicity. An alternating Poisson point process is an example of what is sometimes called a *two-point extension* of a Poisson point process. Characterizing when a finite extension of Bernoulli process is Bernoulli is an important question in Ornstein theory that was considered in a series of papers by Rudolph [24–26].

COROLLARY 8. *An alternating stationary renewal point process with a continuous lifetime distribution is finitarily isomorphic to the same process where the points, additionally, independently receive i.i.d. marks that are uniformly distributed in* $[0, 1]$; *moreover, the isomorphism has a finitary inverse.*

The following theorem relies on the finitary isomorphisms (which also have finitary inverses) constructed by Keane and Smorodinsky [16].

THEOREM 9. *Let $(\Gamma, \mathfrak{X})$ be an independently marked stationary renewal point processes with a continuous lifetime distribution and a skeleton that is a mixing finite state Markov chain with at least two states. There exists a finitary isomorphism between $(\Gamma, \mathfrak{X})$ and a point process with the same distribution as $\Gamma$ that is independently marked with i.i.d. random variables that are uniformly distributed in $[0, 1]$; moreover, the isomorphism has a finitary inverse.*

PROOF. It follows from Keane and Smorodinsky's finitary isomorphism [16] that we may assume $\mathfrak{X}$ is an i.i.d. sequence rather than a general Markov chain. By Lemma 5, together with an elementary argument [31], Lemma 10, we have that the marked point process $(\Gamma, \mathfrak{X})$ is finitarily isomorphic to the point process $\Gamma$ independently marked with i.i.d. random variables that are uniformly distributed in $[0, 1]$, where the inverse is also finitary. □

**4. Reflected Brownian motions.** In this section we recall the definition of reflected Brownian motion and discuss its relation with alternating renewal point processes.

4.1. *Introduction.* Given $h > 0$, we define a Brownian motion $B = (B_s)_{s \in \mathbb{R}}$ *reflected* on the region $[0, h]$ as follows. Let $W = (W_s)_{s \in \mathbb{R}}$ be a Brownian motion where $W_0$ is a real-valued random variable. We set

$$B_s := |W_s - 2nh| \text{ where } n = n(s) \in \mathbb{Z} \text{ is chosen so that} |W_s - 2nh| \leq h.$$

Let

$$p_t(x, y) = \frac{1}{\sqrt{2\pi t}} \exp\left(-\frac{(x-y)^2}{2t}\right)$$

be the transition density of the standard Brownian motion, then the transition density of $B$ is given by

$$q_t(x, y) = \sum_{n \in \mathbb{Z}} [p_t(x, y + 2nh) + p_t(x, -y - 2nh)].$$

For all intervals $(a, b) \subset [0, h]$, we have

$$\lim_{t \to \infty} \int_{-\infty}^{\infty} q_t(x, y) \mathbf{1}_{(a,b)}(y) \, dy = b - a.$$

Hence, it follows that $B$ is stationary if and only if $B_0$ is uniformly distributed on $[0, h]$. Routine considerations yield that $B$ is a Feller process and thus enjoys the strong Markov property; see [13], Chapter 19, for more details. Brownian motion reflected on arbitrary interval $[a, b]$ is defined in a similar way.

4.2. *Point processes within reflected Brownian motions.* Let $h > 0$. Let $B = (B_t)_{t \in \mathbb{R}}$ be a Brownian motion reflected on $[0, h]$. Consider the following random partition of the real line. The real line is partitioned into *upward* intervals of the form $[a, b)$ where $B_a = 0$ and $B_b = h$ with $b = \inf\{t > a : B_t = h\}$ and *downward* intervals of the form $[c, d)$ where $B_c = h$ and $B_d = 0$ with $d = \inf\{t > c : B_t = 0\}$. Our partition is made up by alternating upward and downward intervals that are maximal with respect to inclusion. One such interval can be found by the following finitary procedure.: Let

$$c := \inf\{t \geq 0 : B_t \in \{0, h\}\}.$$

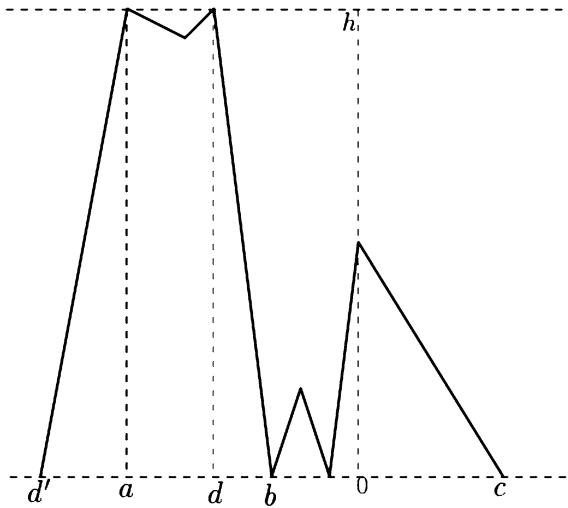

FIG. 1. *The specification of the first interval.*

If $B_c = 0$, then let
$$d := \sup\{t < 0 : B_t = h\} \quad \text{and} \quad b := \inf\{t > d : B_t = 0\};$$
The left endpoint $a$ is defined by another iteration of this procedure, namely,
$$d' := \sup\{t < d : B_t = 0\} \quad \text{and} \quad a := \inf\{t > d' : B_t = h\}.$$
If $B_c = h$, then we define $a$ and $b$ via the same process with the roles of $0$ and $h$ reversed. See Figure 1 for an illustration.

The other intervals are given inductively as follows. Given that we have a finite set of such intervals $\{[b_i, c_i)\}_{i=m}^n$ such that for all $i \in \{m, \ldots, n-1\}$, $b_{i+1} = c_i$ and $[b_i, c_i)$ is an upward or downward interval. Set $b_{n+1} = c_n$ and
$$c_{n+1} := \inf\{t > b_{n+1} : B_t \in \{0, h\} \setminus B_{b_{n+1}}\}.$$
Set also $c_{m-1} = b_m$,
$$d := \sup\{t < b_m : B_t \in \{0, h\} \setminus B_{b_m}\} \quad \text{and} \quad d' := \sup\{t < d : B_t = B_{b_m}\}.$$
Finally, set
$$b_{m-1} := \sup\{t < d' : B_t \in \{0, h\} \setminus B_{b_m}\}.$$

We call the left endpoints of upward intervals the *red* points and the right endpoints of the upward intervals *blue* points. An *upward excursion* on $[a, b)$ is the path $(X_t)_{t \in [a,b)}$ that has the law of Brownian motion reflected on $[0, h]$ conditioned to have an upward interval $[a, b)$. We similarly define *downward excursions*. The following key observation is important in the proofs of our main theorems:

LEMMA 10. *Let $h > 0$. The point process of red and blue points corresponding to a stationary Brownian motion $B$ reflected on the region $[0, h]$ is an alternating renewal point process with lifetime distribution $Y = Y_h$ with the following Laplace transform:*

(1) $$\mathbb{E}(e^{-sY}) = \frac{1}{\cosh(h\sqrt{2s})} \quad \text{for all } s > 0.$$

*Furthermore, given the red and blue points, the joint law of the resulting upward and downward excursions is given by a product measure.*

PROOF. Let $W$ be a Brownian motion and $B$ the corresponding Brownian motion reflected on $[0, h]$ as in Section 4.1. Consider the collection of coloured points $(\mathbf{t}, \mathbf{c}) = (t_j, c_j)_{j \in \mathbb{Z}}$ corresponding to the endpoints of the downward and upward intervals, in order, where $t_0 \leq 0$ is the closest point to the origin from the left and $c_j$ is either red or blue. Note that this collection of points defines a partition of $\mathbb{R}$ to intervals of the form $I_j := [t_j, t_{j+1})$. We refer to the path $B|_{I_j}$ corresponding to the interval $I_j$ as the $j$th excursion. It follows from the strong Markov property for $B$ that, given $(\mathbf{t}, \mathbf{c})$, the excursions $(B|_{I_j})_{j \in \mathbb{Z}}$ are independent.

Clearly, $\mathbf{t}$ is a stationary point process. It remains to show that it is a renewal point process with $Y$ as its lifetime distribution and that $\mathbf{c}$ is independent of $\mathbf{t}$. Let $s \in \mathbb{R}$ be a red or blue point so that $B_s = u \in \{0, h\}$. Let $\tilde{W}$ be a standard Brownian motion started at 0. Then, the stopping times given by

$$\tau(B, s) := \inf\{t - s : t > s, B_t = h - u\}$$

and

$$H := \inf\{t > 0 : \tilde{W}_t \in \{-h, h\}\} \tag{2}$$

have the same distribution, *independent* of the value of $u \in \{0, h\}$. Recall that $T = (T_t)_{t \in \mathbb{R}}$ is the group of translation of $\mathbb{R}$. Since, in addition for all $i \in \mathbb{Z}^+$,

$$t_i = \tau(T_{t_{i-1}}(B), 0) + t_{i-1},$$

it follows that $(t_{i+1} - t_i)_{i \in \mathbb{Z}^+}$ is an i.i.d. sequence whose lifetime distribution is equal to the distribution of $H$ and that $\mathbf{t}$ is independent of $\mathbf{c}$.

By a standard argument [5], Exercise 8.5.1, for all $\theta > 0$, we have that $\exp(\theta W_H - \frac{\theta^2}{2} H)$ is a martingale which implies that

$$\mathbb{E}(e^{-sH}) = \frac{1}{\cosh(h\sqrt{2s})},$$

for all $s > 0$. □

We will refer to the (unmarked) renewal point process in Lemma 10 as a *Brownian excursion point process of parameter $h$*. Thus, the red and blue points are an alternating Brownian excursion point process. We will also need to appeal to the following version of the Borel isomorphism theorem [33], Theorem 3.4.23, for probability spaces:

LEMMA 11. *An upward or downward excursion can be constructed as a bijective function of a random variable that is uniformly distributed on $[0, 1]$.*

PROOF. An upward or downward excursion on the interval $[a, b]$ is a nonatomic probability distribution on the space of a continuous functions with domain $[a, b]$ taking values in $[0, h]$, denoted by $C([a, b], [0, h])$. The space $C([a, b], [0, h])$, endowed with the uniform topology, is a Polish space so that the conclusion follows from the Borel isomorphism theorem. □

PROPOSITION 12. *Let $h > 0$. A stationary Brownian motion reflected on $[0, h]$ is finitarily isomorphic to an alternating Brownian excursion point process of parameter $h$; moreover, the isomorphism has a finitary inverse.*

PROOF. Let $\Gamma$ be an alternating Brownian excursion point process of parameter h. By Lemma 10, $\Gamma$ has the law of the point process of red and blue points corresponding to a Brownian motion reflected on $[0, h]$; thus, all that is missing are the independent downward and upward processes. By Corollary 8, we may assume that the points are independently marked with i.i.d. random variables that are uniformly distributed in [0, 1] and by Lemma 11; there is a bijective function between these uniforms and desired upward and downward excursions. □

PROOF OF THEOREM 4. As in the proof of Proposition 12, a stationary Brownian motion reflected on $[0, h]$ is finitarily isomorphic to an alternating Brownian excursion point process of parameter $h$, where the points are independently marked with i.i.d. random variables that are uniformly distributed in [0, 1]. The Brownian excursion point process gives a partition of $\mathbb{R}$ into upward and downward intervals, and it follows from [32], Lemma 16, that by placing independent Poisson point processes of the same intensity into these intervals we will obtain the desired point process. It is elementary to accomplish this as an explicit function of the uniform markings; see, for example, [32], Proposition 14. □

**5. Periodic Brownian motions.** We define similar machinery to analyze periodic Brownian motions. Let $h > 0$. Let $W = (W_t)_{t \in \mathbb{R}}$ be Brownian motion. We let

$$Z_t := W_t \bmod h$$

so that the points $0$ and $h$ are identified. Then, $Z$ is a *periodic* Brownian motion on $[0, h]$. Note that $Z$ is also a Feller process and thus enjoys the strong Markov property. Consider the following partition of $\mathbb{R}$. Let $n \geq 2$ be an integer. We say $\ell \in [0, h)^n$ is a *vector of levels* if

$$0 = \ell_1 < \ell_2 < \cdots < \ell_n < h.$$

A $(\ell_i, \ell_j)$-*interval* is an interval of the form $[a, b)$ where

$$i - j = \pm 1 \bmod n,$$

$Z_a = \ell_i$ and $Z_b = \ell_j$ with

$$b = \inf\{t > a : B_t = \ell_k \text{ for some } k \text{ with } i - k = \pm 1 \bmod n\}.$$

As with upward and downward intervals, the real line is partitioned into such intervals that are maximal with respect to inclusion. Given vectors of levels $\ell = \{\ell_i\}_{i=1}^n$, we write $\ell_{n+1} = 0$ and $\ell_{-1} = \ell_n$; we find one such interval in the following way: Let

$$c := \inf\{t \geq 0 : B_t \in \{\ell_i\}_{i=1}^n\}.$$

If $B_c = \ell_i$ for $1 \leq i \leq n$, then let

$$d := \sup\{t < 0 : B_t \in \{\ell_{i-1}, \ell_{i+1}\}\} \quad \text{and} \quad b := \inf\{t > d : B_t = \ell_i\}.$$

The left endpoint $a$ is defined by another iteration of this procedure, namely, writing $\ell_j := B_d$; set

$$d' := \sup\{t < d : B_t \in \{\ell_{j-1}, \ell_{j+1}\}\} \quad \text{and} \quad a := \inf\{t > d' : B_t = \ell_j\}.$$

See Figure 2 for an illustration.

The other intervals are given inductively as follows. Given that we have a finite set of such intervals $\{[b_i, c_i)\}_{i=m}^n$ such that for all $k \in \{m, \ldots, n-1\}$, $b_{k+1} = c_i$ and $[b_k, c_k)$ is a $(\ell_{i_k}, \ell_{j_k})$-interval. Then,

$$c_{n+1} := \inf\{t > c_n : B_t \in \{\ell_{j_n-1}, \ell_{j_n+1}\}\}.$$

Set also $c_{m-1} = b_m$,

$$d := \sup\{t < b_m : B_t \in \{\ell_{i_m-1}, \ell_{i_m+1}\}\}, \ell_J := B_d$$

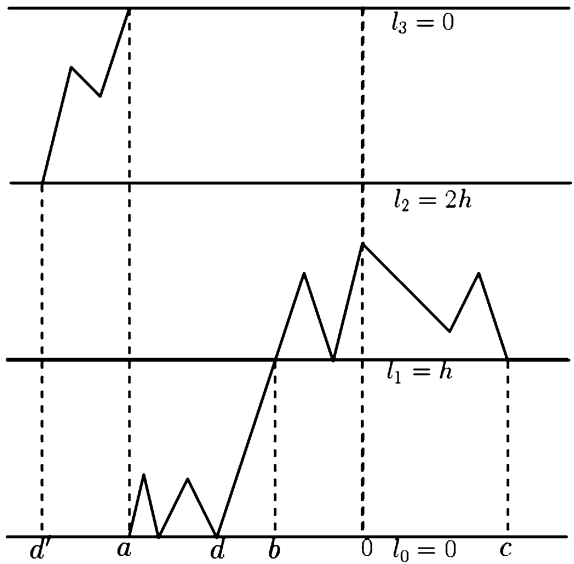

FIG. 2. *The specification of the first interval.*

and
$$d' := \sup\{t < d : B_t \in \{\ell_{J-1}, \ell_{J+1}\}\}.$$

Finally, set
$$b_{m-1} := \sup\{t < d' : B_t = \ell_J\}.$$

An $\ell_i$-*point* is an endpoint $a \in \mathbb{R}$ in the interval of the partition, where $Z_a = \ell_i$; we mark each $\ell_i$ point with the marking $\ell_i$. We call the marked point process of $\ell$-points the *level* point process corresponding to $Z$. An *(i,j)-process* on $[a, b]$ is the process $(Z_t)_{t \in [a,b]}$ that has the law of periodic Brownian motion on $[0, h]$, conditioned to have an $(i, j)$-interval $[a, b]$; we refer to these processes as the *transition* processes. In the special case where $\ell$ is just vector of two points, we refer to a $\ell_1$-point as a *red* point and a $\ell_2$-point as a *blue* point.

LEMMA 13. *Transition processes can be constructed as a bijective function of a random variable that is uniformly distributed on* $[0, 1]$.

PROOF. The proof is identical to the one for Lemma 11. □

LEMMA 14. *Let $h > 0$. Let $Z = (Z_t)_{t \in \mathbb{R}}$ be a stationary periodic Brownian motion on $[0, 2h]$. Let $\ell = (0, h)$ be a vector of two levels. Then, the level point process corresponding to $Z$ is an alternating Brownian excursion point process of parameter $h$. Furthermore, given the red and blue points, the joint law of the resulting transition processes is given by a product measure.*

PROOF. Let $(t_i)_{i \in \mathbb{Z}}$ be the points in the marked point process of $\ell$ points ordered in the usual way. Let $s$ be a red or blue point so that $Z_s = u \in \{0, h\}$, then
$$\tau(Z, s) = \inf\{t - s : t > s, Z_t = h - u\}$$
is distributed as $H$ in (2), regardless of the value of $u \in \{0, h\}$. The rest of the proof is identical to the argument in the proof of Lemma 10 with $Z$ replacing $B$. □

We now state the analogue of Proposition 12 for periodic Brownian motions; the proofs are similar.

PROPOSITION 15. *Let $h > 0$. A stationary periodic Brownian motion on $[0, 2h]$ is finitarily isomorphic to an alternating Brownian excursion point process of parameter $h$; moreover, the isomorphism has a finitary inverse.*

PROOF. Let $\Gamma$ be an alternating Brownian excursion point process with parameter $h$. By Lemma 14, $\Gamma$ has the law of the point process of red and blue points corresponding to a periodic Brownian motion on $[0, 2h]$; thus, all that is missing are the conditionally independent transition processes. By Corollary 8, we may assume that the points also have i.i.d. marks that are uniformly distributed in $[0, 1]$, and, by Lemma 13, there is a bijective function between these uniforms and desired transition processes □

PROOF OF THEOREM 2. Follows immediately from Propositions 12 and 15. □

In order to prove Theorem 3, we will need to consider $n \geq 3$ equally spaced levels.

LEMMA 16. *Let $h > 0$. Let $n \geq 3$ be an integer. Let $Z = (Z_t)_{t \in \mathbb{R}}$ be a stationary periodic Brownian motion on $[0, nh]$. Let $\ell = (h(k-1))_{k=1}^{n}$. Then, the level point process corresponding to $Z$ and $\ell$ is an independently marked Brownian excursion point process with parameter $h$ where the distribution of the markings is an irreducible Markov chain with transition matrix given by*

$$P_{\ell_i, \ell_j} := \begin{cases} \frac{1}{2} & (i - j) = 1 \bmod n, \\ 0 & otherwise, \end{cases}$$

*for all $i, j \in \{1, \ldots, n\}$. Furthermore, conditioned on the level points, the joint law of the resulting transition processes is given by a product measure.*

PROOF. Let $W$ be a Brownian motion and $Z$ the corresponding periodic Brownian motion on $[0, nh]$. Let $(\mathbf{t}, \mathbf{c}) = (t_i, c_i)_{i \in \mathbb{Z}}$ be the level point process corresponding to $Z$, ordered in the usual way, where the markings are given by $c_i \in \{\ell_j\}_{j=1}^{n}$. As in the proof of Lemma 10, this collection of points defines a partition of $\mathbb{R}$ into intervals of the form $I_j := [t_j, t_{j+1})$. We refer to the path $Z|_{I_j}$, corresponding to the interval $I_j$, as the $j$th transition. It follows from the strong Markov property for $Z$ that, given $(\mathbf{t}, \mathbf{c})$, the transitions $(Z|_{I_j})_{j \in \mathbb{Z}}$ are independent.

As before, writing $\ell_0 = h(n-1)$ and $\ell_{n+1} = 0$, then, for every $1 \leq i \leq n$ for a level point $s \in \mathbb{R}$, we have $Z_s = \ell_i$, then

$$\tau(Z, s) := \inf\{t - s : t > s, Z_s \in \{\ell_{i-1}, \ell_{i+1}\}\}$$

has the same law as $H$ as given in (2), regardless of the value of $i \in \{1, \ldots, n\}$. In addition, for all $i \in \mathbb{Z}^+$, we have

$$t_i = \tau(T_{t_{i-1}}(Z), 0) + t_{i-1},$$

so it follows from the strong Markov property for $Z$ that $(t_i - t_{i-1})_{i \in \mathbb{Z}^+}$ is an i.i.d. sequence which is distributed as $H$ and the points $\mathbf{t}$ are independent of the marking sequence $\mathbf{c}$. Since $H$ is as given in (2), the points $\mathbf{t}$ form a Brownian excursion process of parameter $h$.

It remains to show that the marking sequence is a Markov chain with $P$ as its transition matrix. Let $\tilde{W}$ be a standard Brownian motion started at 0. For all $1 \leq i \leq n$ and $k \in \ell$, we have

$$\mathbb{P}(Z_{t_{i+1}} = k | Z_{t_{i-1}} = \ell_i) = \mathbb{P}(\tilde{W}_H = k - \ell_i) \cdot \mathbf{1}_{[k \in \{\ell_{i-1}, \ell_{i+1}\}]} = \frac{1}{2} \mathbf{1}_{[k \in \{l_{i-1}, l_{i+1}\}]};$$

here, the last equality follows from a classical application of the optional stopping theorem; see, for example, [5], Theorem 8.5.3. □

PROOF OF THEOREM 3. Let $h > 0$. We will show that, for all *even* integers $n \geq 4$, the stationary periodic Brownian motion on the interval $[0, nh]$ is finitarily isomorphic to an independently marked alternating Brownian excursion point process with parameter $h$ and that the isomorphism has a finitary inverse. Since any rational can be expressed as the ratio of two *even* integers, the result follows.

As in the proof of Proposition 15, it follows from Lemma 16 and Lemma 13 that a stationary periodic Brownian motion on $[0, nh]$ is finitarily isomorphic (with a finitary inverse) to a Brownian excursion point process with parameter $h$ that is independently coloured with a skeleton that is an irreducible Markov chain and is further independently marked by an i.i.d. sequence of random variables uniformly distributed in $[0, 1]$. Since $n$ is even, the Markov chain has period two and is simply the direct product of a mixing Markov chain and two points. Thus, the stationary periodic Brownian motion is finitarily isomorphic to an alternating Brownian excursion point process with parameter $h$ that is independently coloured with a skeleton that is a mixing Markov chain and is further independently marked by an i.i.d. sequence of random variables uniformly distributed in $[0, 1]$.

It follows readily from Theorem 9 that this coloured marked alternating point process is finitarily isomorphic to an alternating Brownian excursion point process with parameter $h$ that simply independently receives i.i.d. markings that are uniformly distributed in $[0, 1]$. □

REMARK 17. We did not prove that an *unmarked* Brownian excursion point process is finitarily isomorphic to one that independently receives i.i.d. markings that are uniformly distributed in $[0, 1]$; see also Remark 7. The machinery in [9] was developed for Poisson point process, but similar tools may be available for excursion processes.

## 6. Concluding remarks.

6.1. *Explicit isomorphisms*. Although our isomorphisms are finitary, one may say that they are not entirely explicit since we appeal the Borel isomorphism in our proofs of Lemmas 11 and 13.

One advantage of elementary and explicit isomorphisms are that they are easier to understand. Although Keane and Smorodinsky's factors and isomorphisms [14, 15] enjoy both these properties, their constructions are still quite involved; see also [3, 4, 7]. Recall that, armed with entropy, as an isomorphism invariant Kolmogorov was able to distinguish the two Bernoulli shifts $B(\frac{1}{2}, \frac{1}{2})$ and $B(\frac{1}{3}, \frac{1}{3}, \frac{1}{3})$, resolving a decades old problem, where spectral considerations failed; moreover, Sinai proved that the former Bernoulli shift is a factor of the later Bernoulli shift [34]. When recalling this triumph of entropy theory, we are often asked whether there is a simple proof, say in the spirit of Mešalkin [20] or Blum and Hanson [1], of Sinai's factor theorem for this particularly historically relevant example.

6.2. *A family of reflected Brownian motion*. In defining Brownian motion reflected on $[0, h]$, we started with a Brownian motion $W$. Let $\sigma > 0$. More generally, we can start instead with the process $\sigma W$, leading to a family of Brownian motions reflected on $[0, h]$, indexed by the parameter $\sigma > 0$ and we define the analogous related concepts, such as upward and downward paths. The following routine generalization of Lemma 10 assures us that without loss of generality, for the purposes of our theory, we may take $\sigma = 1$.

LEMMA 18. *Let $h > 0$ and $\sigma > 0$. The point process of red and blue points corresponding to a stationary Brownian motion $B$ reflected on the region $[0, h]$ with parameter $\sigma > 0$ is an alternating Brownian excursion point process of parameter $h/\sigma$. Furthermore, given the red and blue points, the joint law of the resulting upward and downward processes is given by a product measure.*

PROOF. We note that, in the proof of Lemma 10, we similarly say that, for every $\theta > 0$, we have that $\exp(\theta \sigma W_t^2 - \frac{(\theta \sigma)^2}{2} t)$ is a martingale and choosing $\theta$ such that $s = \frac{\theta^2 \sigma^2}{2}$ we deduce that

$$\mathbb{E}(e^{-s H_\sigma}) = \frac{1}{\cosh(\frac{h}{\sigma}\sqrt{2s})},$$

where

$$H_\sigma := \inf\{t > 0 : \sigma \tilde{W}_t \in \{-h, h\}\}$$

and $\tilde{W}$ is a standard Brownian motion started at 0. □

PROPOSITION 19. *Let $h > 0$. A stationary Brownian motion reflected on $[0, h]$ with parameter $\sigma > 0$ is finitarily isomorphic to an alternating Brownian excursion point process of parameter $h/\sigma$; moreover, the isomorphism has a finitary inverse.*

PROOF. The proof is, as in the proof of Proposition 12 with Lemma 10, replaced with Lemma 18 and a version of Lemma 11 for upward and downward paths for the reflected Brownian motion with parameter $\sigma > 0$. □

6.3. *More questions.* Related to Question 2 we ask the following:

QUESTION 3. *Do there exist finitary isomorphisms between Brownian excursion point processes and Poisson point processes?*

QUESTION 4. *Find general conditions for renewal point processes to be finitarily isomorphic to Poisson point processes.*

**Acknowledgments.** We thank Benjy Weiss for helpful discussions and encouragement. We also thank Russell Lyons and Richard Bradley for their helpful comments regarding the Proof of Theorem 3. Finally, we thank the referee who provided many good suggestions.

The first author was funded in part by ISF Grant No. 1570/17.

The second author was funded in part by a General Research Fund, provided by the University of Kansas, where he was a faculty member.